\documentclass[11pt]{article}
\usepackage[T1]{fontenc}
\usepackage{amsthm}
\usepackage{longtable}
\usepackage[hidelinks]{hyperref}
\usepackage{amsmath}
\usepackage{amssymb}
\usepackage{a4wide}
\usepackage{lipsum}
\def\Fraisse{Fra\"{\i}ss\' e}
\usepackage{graphicx} 

\usepackage{todonotes}
\def\str#1{\mathbf {#1}}
\def\cont{^\frown}

\def\vb{\mathbf{1}}
\def\Age{\mathrm{Age}}
\def\vc{\mathbf{2}}
\def\eb{\mathbf{11}}
\def\ec{\mathbf{22}}
\def\tb{\mathbf{111}}
\def\Emb{\mathop{\mathrm{Emb}}\nolimits}
\def\tp#1#2{\mathrm{Tp}_{\str{#1}}(#2)}

\newtheorem{theorem}{Theorem}[section]

\theoremstyle{definition}
\newtheorem{definition}[theorem]{Definition}

\theoremstyle{remark}

\title{Counting big Ramsey degrees of the homogeneous and universal $\str{K}_4$-free graph}

\author{
Jan Hubička\thanks{Department of Applied Mathematics (KAM), Charles University, Ma\-lo\-stranské~nám\v estí 25, Praha 1, Czech Republic E-mail: {\tt hubicka@kam.mff.cuni.cz}. Supported by project 25-15571S of  the  Czech  Science Foundation (GA\v CR)}
\and
Matěj Konečný\thanks{Institute of Algebra, TU Dresden, Dresden, Germany. E-mail: {\tt matej.konecny@tu-dresden.de}. Supported by a project that has received funding from the European Union (Project POCOCOP, ERC Synergy Grant 101071674).  Views and opinions expressed are however those of the authors only and do not necessarily reflect those of the European Union or the European Research Council Executive Agency. Neither the European Union nor the granting authority can be held responsible for them.}
\and
Štěpán Vodseďálek\thanks{Charles University, Ma\-lo\-stranské~nám\v estí 25, Praha 1, Czech Republic E-mail: {\tt stepanvodsedalek@gmail.com}. Also supported by project 25-15571S.}
\and
Andy Zucker\thanks{Department of Pure Mathematics, University of Waterloo, 200 University Ave W, Waterloo, ON, Canada. E-mail: {\tt a3zucker@uwaterloo.ca}. Supported by NSERC grants RGPIN-2023-03269 and DGECR-2023-00412.}
}

\makeatletter

\newcommand{\shorttitle}{\@title}

\def\@maketitle{%
  \newpage
  \begin{center}%
  \let \footnote \thanks
    {\small Proceedings of the 13th European Conference on Combinatorics, Graph Theory and Applications\\ EUROCOMB'25\\
    Budapest, August 25 - August 29, 2025
    }
    \vskip 0.5em
    \rule{\linewidth}{0.04cm}
    \vskip 3.5em
    {\LARGE \textbf{\textsc{\@title}} \par}%
    \vskip 1.5em
    {\textbf{\textsc{(Extended abstract)}} \par}
    \vskip 2.5em%
    {\large
      \lineskip .5em%
      \begin{tabular}[t]{c}%
        \@author
      \end{tabular}\par}%
  \end{center}%
  \par
  }
\makeatother

\usepackage{amsthm}

\begin{document}

\thispagestyle{empty}
\maketitle

\begin{abstract}
	Big Ramsey degrees of  \Fraisse{} limits of finitely constrained free amalgamation classes in finite binary languages have been recently fully characterised by Balko, Chodounský, Dobrinen, Hubička, Konečný, Vena, and Zucker.
	A special case of this characterisation is the universal homogeneous $\str{K}_4$-free graph.  We give a self-contained and relatively compact presentation of this case and compute the actual big Ramsey degrees of small graphs.
\end{abstract}


\section{Introduction}
Given graphs $\str{G}$ and $\str{H}$, we denote by $\Emb(\str{G}, \str{H})$ the set of all embeddings $\str{G}\to \str{H}$. If $\str{H}'$ is another graph and $\ell\leq k< \omega$, we write $\str{H}' \longrightarrow (\str{H})^\str{G}_{k,\ell}$ to denote the following statement:
For every colouring $\chi \colon \Emb(\str{G}, \str{H}') \to \{1, \dots, k\}$ with $k$ colours, there exists an embedding $f \colon \str{H} \to \str{H}'$ such that the restriction of $\chi$ to $\Emb(\str{G}, f(\str{H}))$ takes at most $\ell$ distinct values.

For a countably infinite graph $\str{H}$ and a finite induced subgraph $\str{G}$ of $\str{H}$, the \emph{big Ramsey degree of $\str{G}$ in $\str{H}$} is the least number $D\in \omega$ (if it exists) such that $\str{H} \longrightarrow (\str{H})^\str{G}_{k,D}$ for every $k\in \omega$. We say that $\str{H}$ \emph{has finite big Ramsey degrees} if the big Ramsey degree of every finite subgraph $\str{G}$ of $\str{H}$ exists. Big Ramsey degrees of other kinds of structures (orders, hypergraphs, \ldots) are defined in a complete analogy, see recent surveys for details \cite{dobrinen2021ramsey,hubicka2024survey,hubicka2025twenty}.

The concept of big Ramsey degrees, isolated by Kechris, Pestov, and Todorcevic~\cite{Kechris2005},  originated in the study of colourings of subsets of the order of rationals $(\mathbb Q,\leq)$.
In 1969, Laver introduced a rather general proof technique  to obtain upper bounds on big Ramsey degrees of $(\mathbb Q,\leq)$~\cite{todorcevic2010introduction}.
In 1979, Devlin determined the precise big Ramsey degrees proving, somewhat suprisingly, that the big Ramsey degree of a chain with $n$ elements in $(\mathbb Q,\leq)$ is precisely the \emph{$n$-th odd tangent number}: the $(2n-1)$-th derivative of $\tan(x)$ evaluated at 0, the sequence \href{https://oeis.org/A000182}{A000182} in the On-line Encyclopedia of Integer Sequences (OEIS)~\cite{devlin1979,todorcevic2010introduction}.

Graph $\str{H}$ is \emph{homogeneous} if every isomorphism between finite induced subgraphs of $\str{H}$ extends to an automorphism of $\str{H}$. The
\emph{Rado graph} $\str{R}$ is the (up to isomorphism) unique countable homogeneous graph which is \emph{universal}, that is, every countable graph can be embedded to $\str{R}$.
Similarly, for every $k>2$ there exists an (up to isomorphism) unique countable homogeneous $\str{K}_k$-free graph $\str{R}_k$ such that every countable $\str{K}_k$-free graph can be embedded to $\str{R}_k$. We call $\str{R}_k$ the \emph{countable homogeneous $\str{K}_k$-free graph}.
See e.g.~\cite{Hodges1993}.

Laver's proof  can be adapted to the graph $\str{R}$, and in 2006 this was refined by Laflamme, Sauer, and Vuksanovic~\cite{Laflamme2006} to precisely characterise its big Ramsey degrees.
Big Ramsey degrees of cliques and anticliques are again the odd tangent numbers, and
Larson~\cite{larson2008counting} used a Maple program to compute, for a given $n$, the sum of big Ramsey degrees of all graphs with $n$ vertices, yielding a sequence \href{https://oeis.org/A293158}{A293158} in OEIS.

In 2020, Dobrinen developed new techniques to prove finiteness of big Ramsey degrees of $\str{R}_3$~\cite{dobrinen2017universal} (see~\cite{Hubicka2020CS} for a simpler proof) and later of all graphs $\str{R}_k$, $k\geq 3$~\cite{dobrinen2019ramsey}.
Zucker simplified and further generalized Dobrinen's proof to \Fraisse{} limits of finitely constrained free amalgamation classes in finite binary languages~\cite{zucker2020}
and in 2024, Balko, Chodounský, Dobrinen, Hubička, Konečný, Vena, and Zucker gave a precise characterisation~\cite{Balko2021exact}.
In this generality, even the statement of the characterization is very technically challenging and definitions of~\cite{Balko2021exact} need a careful analysis of every specific case they are applied to. The big Ramsey degrees are determined by a number of special trees called \emph{diaries}. To understand them, the reader needs to internalize approximately 21 definitions up to page 22 of~\cite{Balko2021exact}. A short and self-contained description of big Ramsey degrees of $\str{R}_3$ appears in~\cite{Balko2023}.
In this note we give a similar description of diaries of $\str{R}_4$ with the aim to count them.
\section{Diaries of $\str{K}_4$-free graphs}
We first present the definition and then discuss the intuition behind it.
We fix an \emph{alphabet} $\Sigma=\{0,1,2\}$, denote by $\Sigma^*$ the set of all
finite words in the alphabet $\Sigma$,
and by $|w|$ the length of the word $w$. Given $i<|w|$ we denote by $w_i$ the letter of word at index $i$. Indices start by 0.
For $S\subseteq \Sigma^*$, we let $\overline{S}$ be the set $S$ extended by all prefixes of words in $S$.
Given $\ell\geq 0$, we put $\overline{S}_\ell=\{w\in \overline{S}:|w|=\ell\}$.
A word $w\in S$ is a \emph{leaf} of $S$ if there is no $w'\in S$ extending $w$.
Given a word $w$ and a letter $c\in \Sigma$, we denote by $w\cont c$
the word obtained by adding $c$ to the end of $w$. We also set $S\cont c=\{w\cont c:w\in S\}$.

Given distinct $u,v,w\in \Sigma^*$ with $|u|=|v|=|w|=\ell$, we define the following predicates:
\[
	\begin{aligned}
		\vb(u)     & \equiv \exists_{i<\ell}:u_i=1                                &
		\vc(u)     & \equiv \exists_{i<\ell}:u_i=2                                  \\
		\eb(u,v)   & \equiv \exists_{i<\ell}:u_i=v_i=1                            &
		\ec(u,v)   & \equiv \exists_{i<\ell}:u_i=v_i=2                              \\
		\tb(u,v,w) & \equiv \exists_{i<\ell}:u_i=v_i=w_i=1                        &
		u\perp v   & \equiv \neg\vb(u) \hbox{ or } \neg\vb(v)\hbox { or }\ec(u,v)   \\
	\end{aligned}
\]
\begin{definition}[$\str{K}_4$-free diaries]
	\label{def:K4diary}
	A set $S\subseteq \Sigma^*$ is called a \emph{$\str{K}_4$-free-diary} if no member of $S$ extends any other and precisely one of the following seven conditions is satisfied for every $i$ with $0\leq i< \sup_{w\in S}|w|$:
	\begin{enumerate}
		\setlength\itemsep{0em}
		\item \textbf{Splitting (possibly with new $\vb$):}  There is $w\in \overline{S}_i$ such that
		      $
			      \overline{S}_{i+1}=\overline{S}_i{\cont} 0\cup \{w\}\cont 1.
		      $
		\item \textbf{New $\vb$:}  There is $w\in\overline{S}_i$ such that $\neg \vb(w)$ and
		      $
			      \overline{S}_{i+1}=(\overline{S}_i\setminus \{w\})\cont 0\cup \{w\}\cont 1.
		      $
		\item \textbf{New $\vc$:}  There is $w\in\overline{S}_i$ such that $\vb(w)$, $\neg\vc(w)$ and
		      $
			      \overline{S}_{i+1}=(\overline{S}_i\setminus \{w\})\cont 0\cup \{w\}\cont 2.
		      $
		\item \textbf{New $\eb$:} There are distinct words $v,w\in \overline{S}_i$ with $\vb(v)$, $\vb(w)$ and $\neg\eb(v,w)$
		      such that
		      $
			      \overline{S}_{i+1}=(\overline{S}_i\setminus \{v,w\})\cont 0\cup \{v,w\}\cont 1.
		      $
		\item \textbf{New $\ec$:} There are distinct words $v,w\in \overline{S}_i$ with $\vc(v)$, $\vc(w)$, $\eb(v,w)$, $\neg\ec(v,w)$ and $\tb(u,v,w)$ for every $u\in \overline{S}_i$ satisfying $\eb(u,v)$ and $\eb(u,w)$
		      such that
		      $
			      \overline{S}_{i+1}=(\overline{S}_i\setminus \{v,w\})\cont 0\cup \{v,w\}\cont 2.
		      $
		\item \textbf{New $\tb$:} There are distinct words $u,v,w\in \overline{S}_i$ with $\eb(u,v)$, $\eb(u, w)$, $\eb(v, w)$, $\neg\ec(u,v)$, $\neg\ec(u,w)$, $\neg\ec(v,w)$ and $\neg\tb(u,v,w)$
		      such that
		      $
			      \overline{S}_{i+1}=(\overline{S}_i\setminus \{u,v,w\})\cont 0\cup \{u,v,w\}\cont 1.
		      $
		\item \textbf{Leaf:} There is $w\in \overline{S}_i$ with $\vc(w)$ satisfying:
		      \begin{enumerate}
			      \setlength\itemsep{0em}
			      \item \textbf{No new $\eb$:} For every distinct $u,v\in \{z\in \overline{S}_i\setminus\{w\}:z \not \perp w\}$ it holds that $\eb(u,v)$.
			      \item \textbf{No new $\tb$:} For every distinct $u,v,v'\in \{z\in \overline{S}_i\setminus\{w\}:z \not \perp w\}$
			            satisfying $u\not \perp v$, $v\not \perp v'$ and $u\not \perp v'$ it holds that $\tb(u,v,v')$.
			      \item \textbf{No new $\vc$:} For every $u\in \{z\in \overline{S}_i\setminus\{w\}:z \not \perp w\}$ and $v\in S$, $|v|<i$ such that $w_{|v|}=u_{|v|}=1$ it holds that $\vc(u)$.
			      \item \textbf{No new $\ec$:} For every distinct $u,v\in \{z\in \overline{S}_i\setminus\{w\}:z \not \perp w\}$ such that $\tb(u,v,w)$ it holds that $\ec(u,v)$.
			            Moreover for every distinct $u,u'\in \{z\in \overline{S}_i\setminus\{w\}:z \not \perp w\}$ and $v\in S$, $|v|<i$ such that $w_{|v|}=u_{|v|}=u'_{|v|}=1$ it holds that $\ec(u,u')$.

		      \end{enumerate}
		      Moreover:
		      $
			      \overline{S}_{i+1}=\{z\in S_i\setminus\{w\}:z\perp w\}\cont 0 \cup \{z\in S_i\setminus\{w\}:z\not\perp w\}\cont 1.
		      $
	\end{enumerate}
\end{definition}
If $S$ is a $\str{K}_4$-free-diary then by $\str{G}(S)$ we denote the graph with vertex set $S$  with $u,v\in S$, $|u|<|v|$ forming an edge if and only if $v_{|u|}=1$.
Given a $\str{K}_4$-free graph $\str{G}$, we denote by $T(\str{G})$ the set of all $\str{K}_4$-free diaries $S$ for which $\str{G}(S)$ is isomorphic to $\str{G}$.
\begin{theorem}
	For every finite $\str{K}_4$-free graph $\str{G}$, the big Ramsey degree of $\str{G}$ in $\str{R}_4$ equals $|T(\str{G})|\cdot |\mathrm{Aut}(\str{G})|$.
\end{theorem}

\begin{table}[t]
	\centerline{%
		\setlength{\tabcolsep}{0.12em}
		\begin{tabular}{|l|rl|rl|rl|rl|rl|rl|rl|rl|rl|rl|rl|}
			\hline
			\begin{minipage}[c][6mm][t]{0.1mm}%
\end{minipage}
			            & $\str{K}_1$ &                    & $\str{K}_2$ &                   & $\overline{\str{K}}_2$ &                       & $\str{K}_3$     &                     & $\overline{\str{K}}_3$ &                                 & $\str{P}_2$     &                                 & $\overline{\str{P}}_2$ &                                 \\
			\hline
			$\str{R}$   & 1           &                    & $2\cdot 2$  & \cite{Pouzet1996} & $2\cdot 2$             & \cite{Laflamme2006}   & $16\cdot 3!$    & \cite{Laflamme2006} & $16\cdot 3!$           & \cite{Laflamme2006}             & $40\cdot 2$     & \cite{larson2008counting}       & $40\cdot 2$            & \cite{larson2008counting}       \\
			$\str{R}_3$ & 1           & \cite{komjath1986} & $2\cdot 2$  & \cite{sauer1998}  & $5\cdot 2$             & \cite{Balko2021exact} & 0               &                     & $161\cdot 3!$          & \cite{Balko2021exact,Balko2023} & $50\cdot 2$     & \cite{Balko2021exact,Balko2023} & $128 \cdot 2$          & \cite{Balko2021exact,Balko2023} \\
			$\str{R}_4$ & 1           & \cite{Elzahar1989} & $36\cdot 2$ &                   & $23\cdot 2$            &                       & $22658\cdot 3!$ &                     & $197613\cdot 3!$       &                                 & $160488\cdot 2$ &                                 & $267900\cdot 2$        &                                 \\
			\hline
		\end{tabular}
	}
	\label{t1}
	\caption{Big Ramsey degrees of small graphs in $\str{R}$, $\str{R}_3$ and $\str{R}_4$. $\overline{\str{G}}$ denotes the complement of $\str{G}$. $\str{P}_2$ is a path with 2 edges. Big Ramsey degrees are often defined with respect to copies while we use embeddings. The difference between these two values is the size of the automorphism group of the graph. To prevent misunderstandings, we list values with respect to copies explicitly multiplied by the size of the automorphism group.}
\end{table}
\def\ent#1{\parbox{0.28\textwidth}{\small\vskip1mm\begin{minipage}{0.28\textwidth}\begin{flushright} #1\end{flushright}\end{minipage}\vskip1mm}}

\begin{table}[t]
		\begin{longtable}{|c|r|r|r|}
			\hline
			 & $\str{R}$ (\href{https://oeis.org/A000182}{A000182} in OEIS)  & $\str{R}_3$ & $\str{R}_4$  \\
			\hline
			 $\overline{\str{K}}_1$ & \ent{$1\cdot 1!$} & \ent{$1\cdot 1!$} & \ent{$1\cdot 1!$} \\
		$\overline{\str{K}}_2$ & \ent{$2\cdot 2!$} & \ent{$5 \cdot 2!$}  & \ent{$23\cdot 2!$} \\
			 $\overline{\str{K}}_3$ & \ent{$16\cdot 3!$} & \ent{$161 \cdot 3!$} & \ent{$197\,613\cdot 3!$}\\
			 $\overline{\str{K}}_4$ & \ent{$272\cdot 4!$} & \ent{$134\,397 \cdot 4!$}& \ent{\small{$*$ $272\,252\,729\,538\,223\cdot 4!$}}\\
			 $\overline{\str{K}}_5$ & \ent{$7\,936\cdot 5!$} & \ent{$7\,980\,983\,689 \cdot 5!$}& \ent{\tiny{$*$ $43\,391\,315\,736\,159 \hskip4.1mm \allowbreak690\,773\,738\,687\,637\,585\cdot 5!$}}\\
			 $\overline{\str{K}}_6$ & \ent{$353\,792\cdot 6!$} & \ent{$45\,921\,869\,097\,999\,781 \cdot 6!$}& \ent{\tiny{$*$ $1\,075\,426\,511\,374\,671\,039\,522\,386 \hskip4.1mm \allowbreak 376\,330\,779\,194\,191\,609 \hskip4.1mm \allowbreak 662\,344\,240\,057\,102\,999\cdot 6!$}}\\
			 $\overline{\str{K}}_7$ & \ent{$22\,368\,256\cdot 7!$} & \ent{\tiny $35\,268\,888\,847\,472\,944\,795\,910\,097\cdot 7!$}& ?\\
			 $\overline{\str{K}}_8$ & \ent{$1\,903\,757\,312\cdot 8!$} & \ent{\tiny $4\,885\,777\,205\,485\,902\,177 \hskip4.1mm \allowbreak 648\,027\,702\,583\,670\,093\cdot 8!$} & ?\\
			 $\overline{\str{K}}_9$ & \ent{$209\,865\,342\,976\cdot 9!$} & \ent{\tiny $159\,271\,391\,109\,084 \hskip4.1mm \allowbreak 147\,116\,751\,767\,705\,171 \hskip4.1mm \allowbreak 032\,995\,283\,089\,412\,057\cdot 9!$}& ?\\
			 $\overline{\str{K}}_{10}$ & \ent{{$29\,088\,885\,112\,832\cdot 10!$}} & \ent{\tiny $1\,546\,604\,163\,029 \hskip5.3mm \allowbreak 698\,823\,334\,234\,758\,731 \hskip5.3mm \allowbreak 306\,633\,891\,622\,324\,147 \hskip5.3mm \allowbreak 639\,816\,544\,352\,644\,405\cdot 10!$} & ?\\
			\hline
		\end{longtable}
	\label{t1}
	\caption{Big Ramsey degrees of anti-cliques. Values denoted by $*$ have not yet been verified by an independent implementation and should thus be considered preliminary~\cite{Vodsedalek2025bc}.}
\end{table}
A few known values are listed in Table~\ref{t1}.
It is rather surprising to see such a complex structure of diaries arise from three very natural concepts of homogeneity, forbidding a clique, and the big Ramsey degree.
Big Ramsey theorems always fix an enumeration of vertices and the structure arises from the tree of types which we describe now.

A graph is \emph{enumerated} if its vertex set is $\omega$.
Given an enumerated graph $\str{H}$ and $\ell\in \omega$, we call a finite graph $\str{X}$ a \emph{type of $\str{H}$ on level $\ell$} if the vertex set of $\str{X}$ is $\{0,1,\ldots, \ell-1,t\}$ (where $t$ is called the \emph{type vertex}) and the graph created form $\str{X}$ by removing $t$ is an
induced subgraph of $\str{H}$. Types on level $\ell$ thus correspond to one vertex extensions of $\str{H}\restriction_{\{0,1,\ldots,\ell-1\}}$. Graph $\str{X}$ is called a \emph{type of $\str{H}$} if it is a type of $\str{H}$ on level $\ell$ for some $\ell\in \omega$.

We denote by $\mathbb T_\str{H}$ the set of all types of $\str{H}$.  Given $\str{X},\str{X}'\in \mathbb T_\str{H}$ we put $\str{X}\subseteq \str{X}'$ and call \emph{$\str{X}'$ a successor of $\str{X}$} if $\str{X}$ is an induced subgraph of $\str{X}'$. A successor is \emph{immediate} it it differs by only one vertex. Notice that every type has at most two immediate successors and $\mathbb T_\str{H}$ can be viewed as an infinite tree rooted in the unique type on level $0$.
If $\str{X}$ is on level $\ell$ and $\ell'\leq \ell$, we denote by $\str{X}|_{\ell'}$ the unique type $\str{X}'  \in \mathbb T_\str{H}$ on level $\ell$ satisfying $\str{X}'\subseteq \str{X}$.

Given $n\in \omega$ an \emph{$n$-labeled graph $\str{G}$} is a graph we also denote by $\str{G}$ along with a function $\chi_\str{G}$ assigning every vertex $v$ of $\str{G}$ a \emph{label} $\chi_\str{G}(v)\in \{0,1,\ldots, n-1\}$.
To simplify the discussion below, we will additionally require the vertex sets of $n$-labeled graphs to be disjoint from $\omega$.
Given $n\in \omega$, an $n$-labeled graph $\str{G}$, and types $\str{X}_0,\str{X}_1,\ldots, \str{X}_{n-1}$ of $\str{H}$, all on the same level $\ell\in\omega$, we denote by $\str{G}\oplus(\str{X}_0,\str{X}_1,\ldots,\str{X}_{n-1})$ the unique (non $n$-colored) graph
$\str{G}'$ extending $\str{G}$ by vertices $0,1,\ldots,\ell-1$ such that for every vertex $v\in G$ it holds that the subgraph induced by $\str{G}'$ on $\{0,1,\ldots,\ell-1,v\}$ is isomorphic to $\str{X}_{\chi_\str{G}(v)}$ by renaming $t$ to $v$.
Given a tuple $(\str{X}_0,\str{X}_1,\ldots, \str{X}_{n-1})$ of types of a graph $\str{H}$, all on the same level, we denote by $\Age_\str{H}(\str{X}_0,\str{X}_1\allowbreak ,\ldots,\allowbreak  \str{X}_{n-1})$ the set of all finite $n$-labeled graphs $\str{G}$ such that $\str{G}\oplus(\str{X}_0,\str{X}_1,\ldots,\allowbreak \str{X}_{n-1})$ has an embedding to $\str{H}$.

Given an enumerated graph $\str{H}$ and its vertex $v$, we denote by $\tp{H}{v}$ the \emph{type of $v$ in $\str{H}$} created from $\str{H}\restriction \{0,1,\ldots,v\}$ by renaming $v$ to $t$.
Given an enumerated $\str{K}_4$-free graph $\str{H}$, Zucker's theorem can colour only those subgraphs $\str{A}$ of $\str{H}$ which are simultaneously:
\begin{enumerate}
	\item \textbf{Meet-closed}: $\max\{\ell<u: \tp{\str{H}}{u}|_\ell=\tp{\str{H}}{v}|_\ell\}\in A$ for every $u<v\in A$.
	\item \textbf{Closed for age-changes}: For every $u_0,u_1,\ldots, u_{n-1}\in A$ and every $\ell<\min \{u_0,\allowbreak u_1,\allowbreak \ldots u_{n-1}\}$ such that
	      $
		      \Age_\str{H}\allowbreak (\tp{H}{u_0}|_{\ell+1},\allowbreak{}\tp{H}{u_1}|_{\ell+1},\allowbreak{}\ldots,\allowbreak{} \tp{H}{u_{n-1}}|_{\ell+1})
		      \allowbreak \neq
		      \Age_\str{H}\allowbreak (\tp{H}{u_0}|_{\ell},\allowbreak{}\tp{H}{u_1}|_{\ell},\allowbreak{}\ldots,\allowbreak{} \tp{H}{u_{n-1}}|_{\ell})
	      $
	      we have $\ell\in A$.
\end{enumerate}
Given an arbitrary subgraph $\str{A}$ of $\str{H}$, its \emph{closure} is the (unqiue) inclusion minimal subgraph $\str{B}$ of $\str{H}$ which contains $\str{A}$ and satisfies the conditions above.

Designing a diary for $\str{R}_4$ corresponds to finding an enumerated $\str{K}_4$-free graph $\str{H}$ and an embedding $\varphi\colon\str{R}_4\to \str{H}$ which minimizes, for every finite $\str{K}_4$ graph $\str{A}$, the number of order-preserving-isomorphism types of closures of graphs $\varphi(f(\str{A}))$, $f\in \Emb(\str{A},\str{R}_4)$ which then corresponds exactly to its big Ramsey degree.

Methods for minimizing the number of meets were introduced by Devlin, so the main difficulty is to minimize the number of ways age-changes can occur. Let $\str{H}$ be a universal $\str{K}_4$-free graph.
Given a type $\str{X}$ of $\str{H}$ there are three possible sets $\Age_\str{H}(\str{X})$.  If the vertex $t$ is isolated then $\Age_\str{H}(\str{X})$ consists of all finite (1-coloured) $\str{K}_4$-free graphs.  If the neighborhood of $t$ contains no edges then $\Age_\str{H}(\str{X})$ consists of all finite $\str{K}_3$-free graphs. Finally, if the neighborhood of $t$ contains a triangle then $\Age_\str{H}(\str{X})$ contains only graphs with no edges. The second resp. third case corresponds to the predicate $\vb$ resp. $\vc$.

Similarly, given types $\str{X}_0$ and $\str{X}_1$ of $\str{H}$ on the same level $\ell$, $\Age_\str{H}(\str{X}_0,\str{X}_1)$ inherits the structure of $\Age_\str{H}(\str{X}_2)$ on vertices of label $0$ and $\Age_\str{H}(\str{X}_2)$ on vertices of label $1$. There are three options for structures spanning both labels. Either $\Age_\str{H}(\str{X}_0,\str{X}_1)$ contains triangles with both labels, or it contains only edges with both labels (if there exists $i<\ell$ such that $t$ is connected to $i$ in both $\str{X}_0$ and $\str{X}_1$) or it contains no 2-labeled edges (if there exist $i<j$ connected by an edge in $\str{H}$ such that $t$ is connected to both $i$ and $j$ in both $\str{X}_0$ and $\str{X}_1$). Again, the second resp. third case corresponds to the predicate $\eb$ resp. $\ec$. See Example~4.3.5 of~\cite{Balko2021exactLong} for details.

Finally, given types $\str{X}_0$, $\str{X}_1$, and $\str{X}_2$ on level $\ell$, $\Age_\str{H}(\str{X}_0,\str{X}_1,\str{X}_2)$ inherits the structure of ages of each of the pair of types considered. Moreover, it is possible that $\Age_\str{H}(\str{X}_0,\str{X}_1,\str{X}_2)$ contains triangles spanning all three labels. These triangles are blocked if either there exists $i$ such that $t$ is connected to $i$ in all three types (this is captured) by predicate $\tb$), or the age of one of the pairs already forbids the edge spanning the two labels.

If $\str{H}$ is enumerated, every $\str{X}\in \mathbb T_\str{H}$ is determined by its level and the neighborhood of $t$. We can describe this by a word in the alphabet $\Sigma$. Given a word $w$ and $\ell<|w|$, put $i(w,\ell)=\ell+\{j<i:w_j=2\}$. It describes a type $\str{X}_w$ on level $i(w,|w|)$ constructed as follows. If $w_i=0$ then $i(w,\ell)$ and $t$ are not adjacent. If $w_i=1$ then $i(w,\ell)$ is adjacent to $t$. If $w_i=2$ then there are adjacent vertices $i(w,\ell)$ and $i(w,\ell+1)$ (in $\str{H}$ as well as in $\str{X}_w$) both adjacent to $t$.

Given a $\str{K}_4$-free diary $S$, every word $w\in S$ corresponds to a type of some $\str{K}_4$-free graph $\str H$ in which $\str{R}_4$ is embedded.  Leaves correspond to types of vertices of $\str{R}_4$ while non-leaves represent gadgets which are used to reduce the number of closures of subgraphs. These gadgets are of two kinds: a vertex or an edge connected to certain types. Each gadget represents either a meet (splitting) or an age-change, and every change is minimal (so, for example, $\vb$ must happen before $\vc$). Finally, the conditions on leaves signify the fact that ages of every other type with the leaf vertex should already be minimized.

\bibliographystyle{plain}
\bibliography{ramsey}

\begin{thebibliography}{10}

\bibitem{Balko2023}
Martin Balko, David Chodounsk{\' y}, Natasha Dobrinen, Jan Hubi{\v c}ka, Mat{\v
  e}j Kone{\v c}n{\' y}, Llu{\'\i}s Vena, and Andy Zucker.
\newblock Characterisation of the big {R}amsey degrees of the generic partial
  order.
\newblock Submitted, arXiv:2303.10088, 2023.

\bibitem{Balko2021exact}
Martin Balko, David Chodounsk{\' y}, Natasha Dobrinen, Jan Hubi{\v c}ka, Mat{\v
  e}j Kone{\v c}n{\' y}, Llu{\'\i}s Vena, and Andy Zucker.
\newblock Exact big {R}amsey degrees for finitely constrained binary free
  amalgamation classes.
\newblock {\em Journal of the European Mathematical Society}, August 2024.

\bibitem{Balko2021exactLong}
Martin Balko, David Chodounský, Natasha Dobrinen, Jan Hubička, Matěj
  Konečný, Lluis Vena, and Andy Zucker.
\newblock Exact big ramsey degrees for finitely constrained binary free
  amalgamation classes, 2023.
\newblock Manuscript with more discussion and examples than published version,
  arXiv:2110.08409v2.

\bibitem{devlin1979}
Denis Devlin.
\newblock {\em Some partition theorems and ultrafilters on $\omega$}.
\newblock PhD thesis, Dartmouth College, 1979.

\bibitem{dobrinen2017universal}
Natasha Dobrinen.
\newblock The {R}amsey theory of the universal homogeneous triangle-free graph.
\newblock {\em Journal of Mathematical Logic}, 20(02):2050012, 2020.

\bibitem{dobrinen2019ramsey}
Natasha Dobrinen.
\newblock The {R}amsey theory of {H}enson graphs.
\newblock {\em Journal of Mathematical Logic}, 23(01):2250018, 2023.

\bibitem{dobrinen2021ramsey}
Natasha Dobrinen.
\newblock Ramsey theory of homogeneous structures: current trends and open
  problems.
\newblock In {\em International Congress of Mathematicians}, pages 1462--1486.
  EMS Press, 2023.

\bibitem{Elzahar1989}
Mohamed El-Zahar and Norbert~W. Sauer.
\newblock The indivisibility of the homogeneous {$K_n$}-free graphs.
\newblock {\em Journal of Combinatorial Theory, Series B}, 47(2):162--170,
  1989.

\bibitem{Hodges1993}
Wilfrid Hodges.
\newblock {\em Model theory}, volume~42.
\newblock Cambridge University Press, 1993.

\bibitem{Hubicka2020CS}
Jan Hubi{\v{c}}ka.
\newblock Big {R}amsey degrees using parameter spaces.
\newblock {\em Advances in Mathematics}, 2025.
\newblock To appear, arXiv:2009.00967.

\bibitem{hubicka2025twenty}
Jan Hubi{\v{c}}ka and Mat{\v{e}}j Kone{\v{c}}n{\'y}.
\newblock Twenty years of {N}e{\v s}et{\v r}il's classification programme of
  {R}amsey classes.
\newblock {\em arXiv:2501.17293}, 2025.

\bibitem{hubicka2024survey}
Jan Hubi{\v{c}}ka and Andy Zucker.
\newblock A survey on big {R}amsey structures.
\newblock {\em arXiv:2407.17958}, 2025.

\bibitem{Kechris2005}
Alexander~S. Kechris, Vladimir~G. Pestov, and Stevo Todor{\v c}evi{\' c}.
\newblock Fra{\"\i}ss{\'e} limits, {R}amsey theory, and topological dynamics of
  automorphism groups.
\newblock {\em Geometric and Functional Analysis}, 15(1):106--189, 2005.

\bibitem{komjath1986}
P{\'e}ter Komj{\'a}th and Vojtech R{\"o}dl.
\newblock Coloring of universal graphs.
\newblock {\em Graphs and Combinatorics}, 2(1):55--60, 1986.

\bibitem{Laflamme2006}
Claude Laflamme, Norbert~W. Sauer, and Vojkan Vuksanovic.
\newblock Canonical partitions of universal structures.
\newblock {\em Combinatorica}, 26(2):183--205, 2006.

\bibitem{larson2008counting}
Jean~A. Larson.
\newblock Counting canonical partitions in the random graph.
\newblock {\em Combinatorica}, 28(6):659--678, 2008.

\bibitem{Pouzet1996}
Maurice Pouzet and Norbert Sauer.
\newblock Edge partitions of the {R}ado graph.
\newblock {\em Combinatorica}, 16:505--520, 12 1996.

\bibitem{sauer1998}
Norbert~W. Sauer.
\newblock Edge partitions of the countable triangle free homogeneous graph.
\newblock {\em Discrete Mathematics}, 185(1-3):137--181, 1998.

\bibitem{todorcevic2010introduction}
Stevo Todorcevic.
\newblock {\em Introduction to {R}amsey spaces}, volume 174.
\newblock Princeton University Press, 2010.

\bibitem{Vodsedalek2025bc}
Štěpán Vodseďálek.
\newblock Counting big {R}amsey degrees.
\newblock Bachelor's thesis (in preparation), Charles University, 2025.

\bibitem{zucker2020}
Andy Zucker.
\newblock On big {R}amsey degrees for binary free amalgamation classes.
\newblock {\em Advances in Mathematics}, 408:108585, 2022.

\end{thebibliography}

\end{document}